\magnification=1200
\overfullrule=0pt
\centerline {\bf Well-posedness of constrained minimization problems
via saddle-points}\par
\bigskip
\bigskip
\centerline {BIAGIO RICCERI}\par
\bigskip
\bigskip
\noindent
{\it Dedicated to Professor Jean Saint Raymond on his sixtieth birthday, with my greatest
admiration and esteem}\par
\bigskip
\bigskip
\bigskip
\bigskip
Here and in the sequel, $X$ is a Hausdorff topological space, $J, \Phi$ are two
real-valued functions defined in $X$, and $a, b$ are two numbers in $[-\infty,+\infty]$,
with $a<b$. \par
\smallskip
If $a\in {\bf R}$ (resp. $b\in {\bf R}$), we denote by $M_a$ (resp. $M_b$) 
the set of all global minima of the function $J+a\Phi$ (resp. $J+b\Phi$), while
if $a=-\infty$ (resp. $b=+\infty$), $M_a$ (resp. $M_b$) stands for the empty set.
We adopt the conventions $\inf\emptyset=+\infty$, $\sup\emptyset=-\infty$.
\smallskip
We also set
$$\alpha:=\max\left \{ \inf_X \Phi,\sup_{M_b}\Phi\right \}\ ,$$
$$\beta:=\min\left \{ \sup_X \Phi,\inf_{M_a}\Phi\right \}\ .$$
Note that, by Proposition 1 below, one has $\alpha\leq \beta$.\par
\smallskip
A usual, given a function $f:X\to {\bf R}$ and a set $C\subseteq X$, 
  we say that
the problem of minimizing $f$ over $C$ is well-posed if the following two
conditions hold:\par
\smallskip
\noindent
-\hskip 7pt the restriction of $f$ to $C$ has a unique global minimum,
say $\hat x$\ ;\par
\smallskip
\noindent
-\hskip 7pt every sequence $\{x_n\}$ in $C$ such that $\lim_{n\to \infty}f(x_n)=
\inf_Cf$, converges to $\hat x$.\par
\smallskip
A set of the type $\{x\in X: f(x)\leq r\}$ is said to be a sub-level set of
$f$. Clearly, when 
the sub-level sets of $f$ are sequentially compact, the problem of minimizing
$f$ over a sequentially closed set
$C$ is well-posed if and only if $f_{|C}$ has a unique global
minimum.\par
\smallskip
The aim of the present paper is to establish the following result:\par
\medskip
THEOREM 1. - {\it 
Assume that $\alpha<\beta$ and that, for each $\lambda\in
]a,b[$, the function $J+\lambda\Phi$ has sequentially
compact sub-level sets and admits a unique global minimum in
$X$.\par
Then, for each $r\in ]\alpha,\beta[$, 
the problem of minimizing $J$ over $\Phi^{-1}(r)$ is well-posed.\par
Moreover,
if we denote by $\hat x_r$  the unique global minimum
of $J_{|\Phi^{-1}(r)}$ $(r\in ]\alpha,\beta[)$,
the functions $r\to \hat x_r$ and $r\to J(\hat x_r)$ are continuous in
$]\alpha,\beta[$.}
\par
\medskip
Theorem 1 should be regarded as the definitive abstract 
result coming out from the saddle-point method developed in [4], [5], [6], [7], in
specific settings. 
\par
\smallskip
The main tool used to prove Theorem 1 is provided by the following
mini-max result:\par
\medskip
THEOREM 2. -
{\it Let 
$I\subseteq {\bf R}$ be an interval
and $f$ a real-valued function defined in $X\times I$. 
Assume that there exists a number $\rho^*>
\sup_I\inf_X f$, and a point $\hat \lambda\in I$ such that,
for each $\rho\leq\rho^*$, 
 the following conditions hold:\par
\noindent
$(i)$\hskip 5pt the set $\{\lambda\in I : f(x,\lambda)>\rho\}$ is connected
 for all $x\in X$\ ;\par
\noindent
$(ii)$\hskip 5pt
the set $\{x\in X : f(x,\lambda)\leq \rho\}$ is sequentially
closed for all $\lambda\in I$ and sequentially compact
for $\lambda=\hat \lambda$\ ;\par
\noindent
$(iii)$ \hskip 5pt for each compact interval $T\subseteq I$
for which $\sup_T\inf_X f<\rho$, there exists a continuous function
$\varphi:T\to X$ such that $f(\varphi(\lambda),\lambda)<\rho$
for all $\lambda\in T$.\par
Then, one has
$$\sup_{\lambda\in I}\inf_{x\in X}f(x,\lambda)=
\inf_{x\in X}\sup_{\lambda\in I}f(x,\lambda)\ .$$}\par
\smallskip
PROOF. We strictly follow the proof Theorem 2 of [3].
First,
fix a non-decreasing sequence $\{I_n\}$ of compact sub-intervals of $I$,
with $\hat \lambda\in I_1$,
such that $\cup_{n\in {\bf N}}I_{n}=I$. Now, fix $n\in {\bf N}$. We
claim that
$$\sup_{\lambda\in I_n}\inf_{x\in X}f(x,\lambda)=
\inf_{x\in X}\sup_{\lambda\in I_n}f(x,\lambda)\ . \eqno{(1)}$$
Arguing by contradiction, suppose that
$$\sup_{\lambda\in I_n}\inf_{x\in X}f(x,\lambda)<
\inf_{x\in X}\sup_{\lambda\in I_n}f(x,\lambda)\ .$$
Fix $\rho$ satisfying
$$\sup_{\lambda\in I_n}\inf_{x\in X}f(x,\lambda)<\rho<
\min\left \{ \rho^*,\inf_{x\in X}\sup_{\lambda\in I_n}f(x,\lambda)\right \}\ .$$
Set 
$$S=\{(x,\lambda)\in X\times I_n : f(x,\lambda)< \rho\}$$
as well as, for each $\lambda\in I_n$,
$$S^\lambda=\{x\in X: (x,\lambda)\in S\}\ .$$
Since $\sup_{I_n}\inf_Xf<\rho$, one
has $S^\lambda\neq \emptyset$ for all $\lambda\in I_n$. 
 Let $I_n=[a_n,b_n]$.
Put
$$A=\left \{ (x,\lambda)\in S : \lambda<b_n\ ,
\hskip 3pt \sup_{s\in ]\lambda,b_n]}f(x,s)>\rho\right \} $$
and
$$B=\left \{ (x,\lambda)\in S : \lambda>a_n\ ,
\hskip 3pt  \sup_{s\in [a_n,\lambda[}f(x,s)>\rho\right \} \ .$$
Observe that $S^{a_n}\times \{a_n\}\subseteq 
A$ and $S^{b_n}\times \{b_n\}\subseteq B$. Indeed, let $x_1\in S^{a_n}$
and $x_2\in S^{b_n}$. Since
$\rho<\inf_X\sup_{I_n}f$, there are $t, s\in
I_n$ such that $\min\{f(x_1,t),f(x_2,s)\}>\rho$. Since
$\sup\{f(x_1,a_n), f(x_2,b_n)\}<\rho$,
 it follows that
$t>a_n$ and $s<b_n$. Consequently, $(x_1,a_n)\in A$ and
$(x_2,b_n)\in B$. Furthermore, observe that if $(x_0,\lambda_0)\in A$ and if
$\mu\in ]\lambda_0,b_n]$ is such that $f(x_0,\mu)>\rho$, then, in view
of $(ii)$, the set
$$(\{x\in X : f(x,\mu)>\rho\}\times [a_n,\mu[)\cap S$$
is sequentially open in $S$, contains $(x_0,\lambda_0)$ and is contained in $A$.
In other words, $A$ is sequentially open in $S$. Analogously, it is seen that
$B$ is sequentially open in $S$. 
We now prove that $S= A\cup B$. Indeed, let $(x,\lambda)\in
S\setminus A$. We have seen above that
$S^{a_n}\times \{a_n\}\subseteq A$, and so
$\lambda>a_n$. If $\lambda=b_n$, the fact that $(x,\lambda)\in B$ has been likewise proved
above. Suppose $\lambda<b_n$.  Thus, we have 
$\sup_{s\in ]\lambda,b_n]}f(x,s)\leq \rho$. From this, it clearly follows that
$\sup_{s\in [a_n,\lambda[}f(x,s)>\rho$ (note that $f(x,\lambda)<\rho$), 
 and so $(x,\lambda)\in B$.
Furthermore, we have $A\cap B=\emptyset$. Indeed, if $(x_1,\lambda_1)\in A\cap B$,
there would be $t_1,s_1\in I_n$, with $t_1<\lambda_1<s_1$, such that 
$\min\{f(x_1,t_1),f(x_1,s_1)\}>\rho$. By $(i)$, the set
$\{s\in I : f(x_1,s)>\rho\}$ is an interval, and so
 we would have $f(x_1,\lambda_1)>\rho$, against the fact
that $(x_1,\lambda_1)\in S$. Now, in view of $(iii)$, consider
a continuous function $\varphi:I_n\to X$ such that
$$f(\varphi(\lambda),\lambda)<\rho$$
for all $\lambda\in I_n$. Let $h:I_n\to X\times I_n$ be defined by setting
$$h(\lambda)=(\varphi(\lambda),\lambda)$$
for all $\lambda\in I_n$. Since $h$ is continuous, the set
$h(I_n)$ is sequentially connected ([2], Theorem 2.2). But, having in mind that
$h(I_n)\subseteq S$ and that
$h(I_n)$ meets both $A$ and $B$ (since $h(a_n)\in A$ and $h(b_n)\in B$), the
properties of $A, B$ proved above would imply that $h(I_n)$ is sequentially
disconnected, a contradiction. So, $(1)$ holds.
Finally, let us prove the theorem. Again arguing by contradiction,
suppose that
$$\sup_{\lambda\in I}\inf_{x\in X}f(x,\lambda)<
\inf_{x\in X}\sup_{\lambda\in I}f(x,\lambda)\ .$$
Choose $r$ satisfying
$$\sup_{\lambda\in I}\inf_{x\in X}f(x,\lambda)<r<\min
\left \{ \rho^*,\inf_{x\in X}\sup_{\lambda\in I}f(x,\lambda)\right \}\ .$$
 For each $n\in {\bf N}$, put
$$C_n=\left \{ x\in X: \sup_{\lambda\in I_n}f(x,\lambda)\leq r\right \}\ .$$ 
Note that $C_n\neq \emptyset$.
Indeed,
otherwise, we would have
$$r\leq \inf_{x\in X}\sup_{\lambda\in I_n}f(x,\lambda)=
\sup_{\lambda\in I_n}\inf_{x\in X}f(x,\lambda)\leq
\sup_{\lambda\in I}\inf_{x\in X}f(x,\lambda)\ .$$
Consequently, $\{C_n\}$ is a non-increasing
sequence of non-empty sequentially closed subsets
of the sequentially compact set $\{x\in X : f(x,\hat \lambda)\leq \rho^*\}$.
Therefore, one has $\cap_{n\in {\bf N}}C_n\neq \emptyset$.
Let $x^*\in \cap_{n\in {\bf N}}C_n$. Then, one has
$$\sup_{\lambda\in I}f(x^*,\lambda)=
\sup_{n\in {\bf N}}\sup_{\lambda\in I_n}f(x^*,\lambda)
\leq r$$ 
and so
$$\inf_{x\in X}\sup_{\lambda\in I}f(x,\lambda)\leq r\ ,$$
a contradiction. The proof is complete.\hfill $\bigtriangleup$
\medskip
We will also use the following proposition.\par
\medskip
PROPOSITION 1 ([4], Proposition 1). - {\it Let $Y$ be a nonempty set, 
$f, g:Y\to {\bf R}$ two functions, and $\lambda, \mu$ two
real numbers, with $\lambda<\mu$. Let $\hat y_{\lambda}$ 
 be a global minimum of the function $f+\lambda g$ and let
$\hat y_{\mu}$ be a global minimum of the function
$f+\mu g$.\par
Then, one has $$g(\hat y_{\mu})\leq g(\hat y_{\lambda})\ .$$ If either
$\hat y_{\lambda}$
or $\hat y_{\mu}$ is strict and $\hat y_{\lambda}\neq \hat y_{\mu}$, then
$$g(\hat y_{\mu})<g(\hat y_{\lambda})\ .$$}\par
\medskip
{\bf Proof of Theorem 1.} First,
for each $\lambda\in ]a,b[$, denote by $\hat
y_\lambda$ the unique global minimum in $X$ of $J+\lambda\Phi$.
 Let us prove that
the function $\lambda\to \hat y_{\lambda}$ is continuous in $]a,b[$.
To this end, fix $\lambda^*\in ]a,b[$. Let $\{\lambda_n\}$ be any
sequence in $]a,b[$ converging to $\lambda^*$ and let $[c,d]\subset
]a,b[$ be a compact interval containing $\{\lambda_n\}$. Fix $\rho>
\sup_{n\in {\bf N}}\inf_{x\in X}(J(x)+\lambda_n\Phi(x))$. Clearly, we
have
$$\bigcup_{\lambda\in [c,d]}\{x\in X : J(x)+\lambda\Phi(x)\leq \rho\}\subseteq$$
$$\subseteq \{x\in X : J(x)+c\Phi(x)\leq \rho\} \cup \{x\in X : J(x)+d\Phi(x)\leq \rho\}\ .$$
 From this, due to the choice of $\rho$, 
we infer that the sequence $\{\hat y_{\lambda_n}\}$ is contained in the
the set on the right-hand side which is clearly sequentially compact. Hence,
there is a subsequence $\{\hat y_{\lambda_{n_k}}\}$ converging to some
$y^*\in X$. Taking into account that the sequence $\{\Phi(\hat y_{\lambda_{n_k}})\}$
is bounded (by Proposition 1) and that 
the function $J+\lambda^*\Phi$ is sequentially lower semicontinuous,
for each $x\in X$, we then have
$$J(y^*)+\lambda^*\Phi(y^*)\leq
\liminf_{k\to \infty}(J(\hat y_{\lambda_{n_k}})+\lambda^*\Phi(\hat y_{\lambda_{n_k}}))=$$
$$=\liminf_{k\to \infty}(J(\hat y_{\lambda_{n_k}})+\lambda_{n_k}\Phi(\hat y_{\lambda_{n_k}})+
(\lambda^*-\lambda_{n_k})\Phi(\hat y_{\lambda_{n_k}}))=$$
$$=\liminf_{k\to \infty}(J(\hat y_{\lambda_{n_k}})+\lambda_{n_k}\Phi(\hat y_{\lambda_{n_k}}))\leq
\lim_{k\to \infty}(J(x)+\lambda_{n_k}\Phi(x))=J(x)+\lambda^*\Phi(x)\ .$$
Hence $y^*$ is the global minimum of $J+\lambda^*\Phi$, that is $y^*=\hat y_{\lambda^*}$, which
shows the continuity of $\lambda\to \hat y_{\lambda}$ at $\lambda^*$. 
Now, fix $r\in ]\alpha,\beta[$
 and consider the function
$f:X\times {\bf R}\to {\bf R}$ defined by
$$f(x,\lambda)=J(x)+\lambda(\Phi(x)-r)$$
for all $(x,\lambda)\in X\times {\bf R}$. Clearly,
the the restriction of the function $f$ to $X\times
]a,b[$
satisfies all the assumptions of Theorem 1. In particular,
$(iii)$ is satisfied taking $\varphi(\lambda)=\hat y_{\lambda}$.
Consequently, we have
$$\sup_{\lambda\in ]a,b[}\inf_{x\in X}(J(x)+\lambda
(\Phi(x)-r))=
\inf_{x\in X}\sup_{\lambda\in ]a,b[} (J(x)+\lambda(\Phi(x)-r))
\ .\eqno{(2)}$$
Note that
$$\sup_{\lambda\in ]a,b[}\inf_{x\in X}f(x,\lambda)\leq
\sup_{\lambda\in [a,b]\cap {\bf R}}\inf_{x\in X}f(x,\lambda)\leq$$
$$\leq\inf_{x\in X}\sup_{\lambda\in [a,b]\cap {\bf R}}
f(x,\lambda)=
\inf_{x\in X}\sup_{\lambda\in ]a,b[}f(x,\lambda)$$
and so from $(2)$ it follows
$$\sup_{\lambda\in [a,b]\cap {\bf R}}\inf_{x\in X}(J(x)+\lambda
(\Phi(x)-r))=
\inf_{x\in X}\sup_{\lambda\in [a,b]\cap {\bf R}} (J(x)+\lambda(\Phi(x)-r))\ .\eqno{(3)}$$
Now, observe that the function $\inf_{x\in X}f(x,\cdot)$ is
upper semicontinuous in $[a,b]\cap {\bf R}$ and that
$$\lim_{\lambda\to +\infty}\inf_{x\in X}f(x,\lambda)=-\infty$$
if $b=+\infty$ (since $r>\inf_X \Phi$), and
$$\lim_{\lambda\to -\infty}\inf_{x\in X}f(x,\lambda)=-\infty$$
if $a=-\infty$ (since $r<\sup_X \Phi)$.
 From this, it clearly follows that
 there exists $\hat \lambda_r\in [a,b]\cap {\bf R}$ such that
$$\inf_{x\in X}f(x,\hat \lambda_r)=\sup_{\lambda\in [a,b]\cap {\bf R}}\inf_{x\in X}f(x,\hat \lambda_r)\ .$$
Since
$$\sup_{\lambda\in [a,b]\cap {\bf R}}f(x,\lambda)=\sup_{\lambda\in ]a,b[}f(x,\lambda)$$
for all $x\in X$, the sub-level sets of the function
$\sup_{\lambda\in [a,b]\cap {\bf R}}f(\cdot,\lambda)$ are sequentially compact.
 Hence, there exists $\hat x_r\in
X$ such that
$$\sup_{\lambda\in [a,b]\cap {\bf R}}f(\hat x_r,\lambda)=
\inf_{x\in X}\sup_{\lambda\in [a,b]\cap {\bf R}}f(x,\lambda)\ .$$
Then, thanks to $(3)$, $(\hat x_r,\hat \lambda_r)$ is a saddle-point of $f$, that is
$$J(\hat x_r)+\hat \lambda_r (\Phi(\hat x_r)-r)=\inf_{x\in X}(J(x)+\hat \lambda_r(\Phi(x)-r))=
J(\hat x_r)+\sup_{\lambda\in [a,b]\cap {\bf R}}\lambda(\Phi(\hat x_r)-r)\ .\eqno{(4)}$$
First of all, from $(4)$ it follows that $\hat x_r$ is a global minimum of
$J+\hat \lambda_r\Phi$.
We now show that $\Phi(\hat x_r)=r$. We
distinguish four cases.\par
\noindent
-\hskip 7pt $a=-\infty$ and $b=\infty$.  In this case, the equality $\Phi(\hat x_r)=
r$ follows from the fact that
$\sup_{\lambda\in {\bf R}}\lambda(\Phi(\hat x_r)-r)$ is finite.\par
\noindent
-\hskip 7pt $a>-\infty$ and $b=+\infty$. In this case, the finiteness of
$\sup_{\lambda\in [a,+\infty[}\lambda(\Phi(\hat x_r)-r)$ implies that $\Phi(\hat x_r)\leq
r$. But, if $\Phi(\hat x_r)<r$, from $(4)$, we would infer that
$\hat \lambda_r=a$ and so $\hat x_r\in M_a$. This would imply
 $\inf_{M_a}\Phi<r$, contrary to the choice of $r$.\par
\noindent
-\hskip 7pt $a=-\infty$ and $b<+\infty$. In this case, the finiteness of
$\sup_{\lambda\in ]-\infty,b]}\lambda(\Phi(\hat x_r)-r)$ implies that
$\Phi(\hat x_r)\geq r$. But, if $\Phi(\hat x_r)>r$, from $(4)$ again, we would infer
 $\hat \lambda_r=b$, and so $\hat x_r\in M_b$. Therefore, $\sup_{M_b}\Phi>r$, contrary to
the choice of $r$.\par
\noindent
-\hskip 7pt $-\infty<a$ and $b<+\infty$. In this case, if $\Phi(\hat x_r)\neq r$, as we
have just seen, we would have either $\inf_{M_a}\Phi<r$ or $\sup_{M_b}\Phi>r$, contrary
to the choice of $r$.\par
Having proved that $\Phi(\hat x_r)=r$, we also get that $\hat \lambda_r\in ]a,b[$.
Indeed, if $\hat \lambda_r\in \{a,b\}$, we would have either
$\hat x_r\in M_a$ or $\hat x_r\in M_b$ and so either $\inf_{M_a}\Phi
\leq r$ or $\sup_{M_b}\Phi\geq r$, contrary to the choice of $r$. From $(4)$ once again,
we furthermore infer that any global minimum of
$J_{|\Phi^{-1}(r)}$ (and $\hat x_r$ is so)
is a global minimum of $J+\hat \lambda_r\Phi$
in $X$. But, since $\hat \lambda_r\in ]a,b[$, $J+\hat \lambda_r\Phi$ has exactly one
global minimum in $X$ which, therefore, coincides with $\hat x_r$. Since
the sub-level sets of $J+\hat \lambda_r\Phi$ are sequentially compact, we then
conclude that any minimizing sequence in $X$ for $J+\hat \lambda_r\Phi$ converges
to $\hat x_r$. But any minimizing sequence in $\Phi^{-1}(r)$ for
$J$ is a minimizing sequence for
$J+\hat \lambda_r\Phi$, and so it converges to $\hat x_r$.
Consequently, the problem of minimizing $J$ over $\Phi^{-1}(r)$
is well-posed, as claimed. \par
Now, let us prove the other assertions made in thesis.
By Proposition 1, it clearly follows that the function $\lambda\to
\Phi(\hat y_{\lambda})$ is non-increasing in $]a,b[$ and that its
range is contained in $[\alpha,\beta]$. On the other hand, by the
first assertion of the thesis, this range contains $]\alpha,\beta[$.
Of course, from this it follows that the function $\lambda\to
\Phi(\hat y_{\lambda})$ is continuous in $]a,b[$. Now, observe
that the function $\lambda\to \inf_{x\in X}(J(x)+\lambda\Phi(x))$ is
concave and hence continuous in $]a,b[$. This, in particular, implies
that the function $\lambda\to J(\hat y_{\lambda})$ is continuous in
$]a,b[$.
Now, for each
$r\in ]\alpha,\beta[$, put
$$\Lambda_r=\{\lambda\in ]a,b[ : \Phi(\hat y_{\lambda})=r\}\ .$$
Let us prove that the multifunction $r\to \Lambda_r$ is upper semicontinuous in $]\alpha,\beta[$.
Of course, it is enough to show that the restriction of the multifunction to any bounded open
sub-interval of $]\alpha,\beta[$ is upper semicontinuous. So, let $s, t\in ]\alpha,\beta[$, with
$s<t$. Let $\mu, \nu\in ]a,b[$ be such that $\Phi(\hat y_{\mu})=t$, $\Phi(\hat y_{\nu})=s$. By
Proposition 1, we have 
$$\bigcup_{r\in ]s,t[}\Lambda_r\subseteq [\mu,\nu]\ .$$
Then, to show that the restriction of multifunction $r\to \Lambda_r$ to $]s,t[$
is upper semicontinuous, it is enough to prove that its graph is closed
in $]s,t[\times [\mu,\nu]$ ([1], Theorem 7.1.16).
But, this latter fact follows immediately from the continuity of the function
$\lambda\to \Phi(\hat y_{\lambda})$. At this point, we observe that, for
each $r\in ]\alpha,\beta[$, the function $\lambda\to \hat y_{\lambda}$ is
constant in $\Lambda_r$. Indeed, let $\lambda,\mu\in \Lambda_r$ with
$\lambda\neq \mu$. If it was $\hat y_{\lambda}\neq \hat y_{\mu}$, by
Proposition 1 it would follow
$$r=\Phi(\hat y_{\lambda})\neq \Phi(\hat y_{\mu})=r\ ,$$
an absurd. Hence, the function $r\to \hat x_r$, as composition of
the upper semicontinuous multifunction $r\to \Lambda_r$ and the continuous
function $\lambda\to \hat y_{\lambda}$, is continuous. Analogously,
the continuity of the function $r\to J(\hat x_r)$ follows observing that it
is the composition of $r\to \Lambda_r$ and the continuous function
$\lambda\to J(\hat y_{\lambda})$. The proof is complete.
\hfill $\bigtriangleup$\par
\medskip
REMARK 1. - We want to point out that, under the assumptions of Theorem 1, we have
actually proved that, for each $r\in ]\alpha,\beta[$, there exists $\hat \lambda_r\in ]a,b[$
such that the unique global minimum of $J+\hat \lambda_r\Phi$ belongs to $\Phi^{-1}(r)$.\par
\medskip
When $a\geq 0$, we can obtain a conclusion dual to
that of Theorem 1,
 under the same key assumption.\par
\medskip
THEOREM 3. - {\it  Let $a\geq 0$.
 Assume that, for each $\lambda\in
]a,b[$, the function $J+\lambda\Phi$ has sequentially
compact sub-level sets and admits a unique global minimum in
$X$.\par
Set
$$\gamma:=
\max\left \{ \inf_X J,\sup_{\hat M_a}J\right \}\ ,$$
$$\delta:=\min\left \{ \sup_X J,
\inf_{\hat M_b}J\right \}\ ,$$
where
$$\hat M_a=\cases {M_a & if $a>0$\cr & \cr \emptyset & if $a=0$\ ,\cr}$$
$$\hat M_b=\cases {M_b & if $b<+\infty$\cr & \cr \inf_X\Phi & if $b=+\infty$\ .\cr}$$
Assume that $\gamma<\delta$.\par
Then, for each $r\in ]\gamma,\delta[$, 
the problem of minimizing $\Phi$ over $J^{-1}(r)$ is well-posed.\par
Moreover,
if we denote by $\tilde x_r$  the unique global minimum
of $\Phi_{|J^{-1}(r)}$ $(r\in ]\gamma,\delta[)$,
the functions $r\to \tilde x_r$ and $r\to \Phi(\tilde x_r)$ are continuous in
$]\gamma,\delta[$.}\par
\smallskip
PROOF. Let $\mu\in ]b^{-1},a^{-1}[$. Then, since $\mu^{-1}\in ]a,b[$
and
$$\Phi+\mu J=\mu(J+\mu^{-1}\Phi)\ ,$$
we clearly have that the function $J+\mu\Phi$ has sequentially compact
sub-level sets and admits a unique global minimum. At this point, the conclusion
follows applying Theorem 1 with the roles of $J$ an $\Phi$ interchanged.
\hfill $\bigtriangleup$ \par
\medskip
We now state the version of Theorem 1 obtained in the setting of a
reflexive Banach space endowed with the weak topology.\par
\medskip
THEOREM 4. - {\it Let $X$ be a sequentially weakly closed set in a reflexive
real Banach space.
Assume that $\alpha<\beta$ and that, for each $\lambda\in
]a,b[$, the function $J+\lambda\Phi$ is sequentially weakly lower
semicontinuous, has bounded sub-level sets and has a unique global minimum in
$X$.  \par
Then, for each $r\in ]\alpha,\beta[$, 
the problem of minimizing $J$ over $\Phi^{-1}(r)$ is well-posed in the
weak topology. \par
Moreover,
if we denote by $\hat x_r$  the unique global minimum
of $J_{|\Phi^{-1}(r)}$ $(r\in ]\alpha,\beta[)$,
the functions $r\to \hat x_r$ and $r\to J(\hat x_r)$ are continuous in
$]\alpha,\beta[$, the first one in the weak topology.}
\smallskip
PROOF. Our assumptions clearly imply that,
for each $\lambda\in ]a,b[$, the sub-level sets of $J+\lambda\Phi$
are sequentially weakly compact, by the Eberlein-
Smulyan theorem. Hence, considering $X$ with the relative weak topology,
we are allowed to apply Theorem 1, from which the conclusion directly
follows.\hfill $\bigtriangleup$\par
\medskip
Analogously, from Theorem 3 we get\par
\medskip
THEOREM 5. - {\it Let $a\geq 0$ and let
 $X$ be a sequentially weakly closed set in a reflexive
real Banach space.
Assume that, for each $\lambda\in
]a,b[$, the function $J+\lambda\Phi$ is sequentially weakly lower
semicontinuous, has bounded sub-level sets and has a unique global minimum in
$X$. Assume also that $\gamma<\delta$, where $\gamma, \delta$ are defined as
in Theorem 3.\par
Then, for each $r\in ]\gamma,\delta[$, 
the problem of minimizing $\Phi$ over $J^{-1}(r)$ is well-posed in the
weak topology.
 \par
Moreover,
if we denote by $\tilde x_r$  the unique global minimum
of $\Phi_{|J^{-1}(r)}$ $(r\in ]\gamma,\delta[)$,
the functions $r\to \tilde x_r$ and $r\to \Phi(\tilde x_r)$ are continuous in
$]\gamma,\delta[$, the first one in the weak topology.}\par
\medskip
Finally, it is worth noticing that Theorem 1 also offers 
the perspective
of a novel way of seeing whether a given function possesses a global
minimum. Let us formalize this using Remark 1.\par
\medskip
THEOREM 6. - {\it
Assume that $b>0$ and that,
 for each $\lambda\in ]0,b[$,
the function $J+\lambda\Phi$ has sequentially compact sub-level sets and
 admits a unique global minimum, say
$\hat y_{\lambda}$. Assume also that
$$\lim_{\lambda\to 0^+}\Phi(\hat y_{\lambda})<\sup_X \Phi\ .\eqno{(5)}$$
Then, one has
$$ \lim_{\lambda\to 0^+}\Phi(\hat y_{\lambda})=\inf_M \Phi\ ,$$
where $M$ is the set of all global minima of $J$ in $X$.}\par
\smallskip
PROOF. We already know that the function $\lambda\to \Phi(\hat y_{\lambda})$
is non-increasing in $]a,b[$ and that its range is contained in
$[\alpha,\beta]$. We claim that
$$\beta=\lim_{\lambda\to 0^+}\Phi(\hat y_{\lambda})\ .$$
Assume the contrary.
 Let us apply Theorem 1, with
$a=0$ (so, $M_0=M$), using the conclusion pointed out in Remark 1. Choose $r$
satisfying
$$\lim_{\lambda\to 0^+}\Phi(\hat y_{\lambda})<r<\beta\ .$$
Then, (since also $\alpha<r$) it would exist $\hat \lambda_r\in ]0,b[$ such that
$\Phi(\hat y_{\hat\lambda_r})=r$,
 contrary to the choice of $r$.
 At this point, the conclusion follows directly from
 $(5)$. 
\hfill $\bigtriangleup$\par
\vfill\eject
\centerline {\bf References}\par
\bigskip
\bigskip
\noindent
[1]\hskip 5pt A. FEDELI and A. LE DONNE, {\it On good connected preimages}, Topology
Appl., {\bf 125} (2002), 489-496.\par
\smallskip
\noindent
[2]\hskip 5pt E. KLEIN and A. C. THOMPSON, {\it Theory of Correspondences}, John Wiley
$\&$ Sons, 1984.\par
\smallskip
\noindent
[3]\hskip 5pt B. RICCERI, {\it Minimax theorems for limits of parametrized functions
having at most one local minimum lying in a certain set}, Topology Appl., 
{\bf 153} (2006), 3308-3312.\par
\smallskip
\noindent
[4]\hskip 5pt B. RICCERI, {\it Uniqueness properties of functionals with
Lipschitzian derivative}, Port. Math. (N.S.), {\bf 63} (2006), 393-400.\par
\smallskip
\noindent
[5]\hskip 5pt B. RICCERI, {\it On the existence and uniqueness of minima and maxima
on spheres of the integral functional of the calculus of variations}, J. Math. 
Anal. Appl., {\bf 324} (2006), 1282-1287.\par
\smallskip
\noindent
[6]\hskip 5pt B. RICCERI, {\it On the well-posedness of optimization problems on
spheres in $H^1_0(0,1)$}, J. Nonlinear Convex Anal., {\bf 7} (2006), 525-528.\par
\smallskip
\noindent
[7].\hskip 5pt B. RICCERI, {\it The problem of minimizing locally a $C^2$ functional
around non-critical points is well-posed}, Proc. Amer. Math. Soc., {\bf 135}
(2007), 2187-2191.\par
\bigskip
\bigskip
Department of Mathematics\par 
University of Catania\par
Viale A. Doria 6\par
95125 Catania\par
Italy\par
\smallskip
{\it e-mail address}: ricceri@dmi.unict.it
\bye